\documentclass[smallextended,referee,envcountsect]{svjour3} 

\begin{document}

\title{A proximal point algorithm revisited and extended}

\subtitle{}

\author{Gheorghe Moro\c{s}anu}

\institute{
Author's Address: \\
Central European University \\ 
Department of Mathematics and its Applications \\
Nador u. 9 \\
1051 Budapest, Hungary\\
morosanug@ceu.edu}

\date{Received: date / Accepted: date}

\maketitle

\begin{abstract}
This Note is inspired by the recent paper by Djafary Rouhani and Moradi [J. Optim. Theory Appl. 172 (2017) 222-235], where a proximal point algorithm 
proposed by Boikanyo and Moro\c{s}anu [Optim. Lett. 7 (2013) 415-420] is discussed. We start with a brief history of the subject and then propose and analyse 
the following more general algorithm for approximating the zeroes of a maximal monotone operator $A$ in real Hilbert space $H$
$$
x_{n+1}=(I+\beta_nA)^{-1}(u_n + \alpha_n(x_n+e_n)), \ \ n\ge 0\, ,  
$$
where $x_0\in H$ is a given starting point, 
$u_n \rightarrow u$ is a given sequence in $H$, 
${R} \ni \alpha_n \rightarrow 0$, and $(e_n)$ 
is the error sequence satisfying 
$\alpha_ne_n\rightarrow 0$. Besides the main result on the strong convergence of $(x_n)$, we discuss some particular 
cases, including the approximation of minimizers of convex functionals, explain how to use our algorithm in practice, and present some simulations to 
illustrate the applicability of our algorithm. 
\end{abstract}


\section{Introduction}

Let $H$ be a real Hilbert space with scalar product $(\cdot , \cdot )$ and norm $\Vert \cdot \Vert$. An operator $A:D(A)\subset H \to H$ (possibly set-valued) 
is said to be monotone if its graph $G(A)=\{ [x,y]\in D(A)\times H; \, y\in Ax\}$ is a monotone subset of $H\times  H$, i.e., 
$$
(x_1-x_2,y_1-y_2)\ge 0 \ \ \forall [x_1,y_1], \, [x_2,y_2]\in G(A)\, . 
$$
If in addition $G(A)$ is not properly contained in the graph of any other monotone operator in $H$, then $A$ is called maximal monotone. It is well-known 
that a monotone operator $A$ is maximal monotone if and only if the range of $I+\lambda A$ is all of $H$ for all $\lambda >0$ (equivalently for 
some $\lambda >0$). In this case the so-called resolvent operator $J_{\lambda}=(I+\lambda A)^{-1}$ is everywhere defined, single-valued and nonexpansive 
(i.e., Lipschitz with constant $L=1$). If $\phi:H\to (-\infty, +\infty]$ is a proper (i.e., not identically $+\infty$), convex, lower semicontiunous 
function then the subdifferential operator defined by 
$$
\partial \phi (x)=\{ y\in H; \, \phi (x)-\phi (v)\le (y,x-v) \ \forall v\in D(\phi)\} 
$$
is maximal monotone. For more information on monotone operators and convex functions see \cite{brezis} and \cite{morosanu}.

We are interested in solving the problem
\begin{equation}\label{eq}
\mbox{Find } x\in D(A)  \ \ \mbox{ such that } 0\in Ax \, .
\end{equation}
Denote by $F$ the solution set of (\ref{eq}), i.e., $F=A^{-1}0$. One of the most important iterative methods for finding approximate 
solutions of (\ref{eq}) is the proximal point algorithm (PPA) which was introduced by Martinet \cite{martinet} for  a particular case of $A$ 
and then extended by Rockafellar \cite{rockafellar} to a general maximal monotone operator. For each $x_0\in H$ the PPA generates the sequence 
$(x_n)$ as follows 
\begin{equation}\label{PPA}
x_{n+1}=J_{\beta_n}(x_n + e_n), \ \ n\ge 0\, , 
\end{equation}
where $\beta_n \in (0, \infty)$ for all $n\ge 0$ and $(e_n)$ is the sequence of computational errors. Unfortunately (under the suitable conditions $\lim \inf \beta_n >0$, $\sum_{n=0}^{\infty}\Vert e_n\Vert <\infty$) 
$(x_n)$ converges in general only weakly (to points of F), even in the particular case when $A$ is a subdifferential operator (see \cite{guler}). Subsequently 
much work has been dedicated towards modifying the PPA to obtain algorithms that generate strongly convergent sequences. Recall that, inspired by 
Lehdili and Moudafi's prox-Tikhonov method (see \cite{lehmud}), Xu \cite{xu} considered the following iterative scheme 
\begin{equation}\label{xualg}
x_{n+1}=J_{\beta_n}(\lambda_nu+(1-\lambda_n)x_n + e_n), \ \ n\ge 0\, , 
\end{equation}
where $\lambda_n\in (0,1), \ \beta_n\in (0,\infty) \ \ \forall n\ge 0$, $\lambda_n\rightarrow 0$, 
$\sum_{n=1}^{\infty}\lambda_n=\infty$, and proved that, under some additional conditions, $x_n$ converges strongly to 
$P_Fu$, the 
metric projection of $u$ onto $F$ (which was assumed to be nonempty). The best convergence result on (\ref{xualg}) has 
been reported later by Wang and Cui \cite{wangcui}. Specifically, they 
proved that $(x_n)$ generated by (\ref{xualg}) converges strongly to $P_Fu$ under the following conditions: 
$F\neq \emptyset$, $\lambda_n\in (0,1)$, $\beta_n\in (0,\infty) \ \forall n\ge 0$, $\lim \inf \beta_n >0, \ 
\lambda_n \rightarrow 0, \ 
\sum_{n=0}^{\infty} \lambda_n = \infty$, and either $\sum_{n=0}^{\infty}\Vert e_n \Vert < \infty$ or $\lim \Vert e_n \Vert /\lambda_n = 0$. In fact, under 
these conditions, (\ref{xualg}) is equivalent with
\begin{equation}\label{xualg2}
x_{n+1}=J_{\beta_n}(\lambda_nu+(1-\lambda_n)(x_n + e_n)), \ \ n\ge 0\, .  
\end{equation}
In \cite{bm} a strong convergence result for $(x_n)$ generated by (\ref{xualg2}) was reported in the case of the 
alternative framework: $F\neq \emptyset$, 
$\lambda_n \in (0,1), \ \beta_n \in (0, \infty) \ \ \forall n\ge 0$, $\lambda_n \rightarrow 1$, $\beta_n \rightarrow \infty$, and $(e_n)$ bounded. The same 
framework is reconsidered in a recent paper by Djafari Rouhani and Moradi \cite{rouhanimoradi}. They use the 
condition $(\lambda_n -1)e_n \rightarrow 0$ (instead of the boundedness of $(e_n)$). In fact this condition is 
also easily visible from our approach in \cite{bm}.  

\bigskip

The main observation leading to this Note is that: while the convex combination $\lambda_nu+(1-\lambda_n)(x_n + e_n)$ in (\ref{xualg2}) is relevant when 
$\lambda_n \rightarrow 0$, it is not the case if $\lambda_n \rightarrow 1$. Indeed, we can consider the following more general algorithm
\begin{equation}\label{generalalg}
x_{n+1}=J_{\beta_n}(u_n+\alpha_n(x_n+e_n)), \ \ n\ge 0 \, , 
\end{equation}
where 

$(H)$ ~~~ $A:D(A)\subset H\to H$ is a maximal monotone operator with $A^{-1}0=:F \neq \emptyset$; $\beta_n \in (0,\infty)$, $\alpha_n\in {R}$ for 
all $n\ge 0$, $\beta_n \rightarrow \infty$, $\alpha_n \rightarrow 0$; $\alpha_n e_n \rightarrow 0$; $u_n \rightarrow u$. 

\bigskip

Our main result (Theorem \ref{main}) states that under hypotheses $(H)$, for every $x_0\in H$, the sequence $(x_n)$ generated by (\ref{generalalg})  
converges strongly to $P_Fu$. By chosing $u_n=\lambda_n u$ and $\alpha_n = 1-\lambda_n, \ n\ge 0$ with $\lambda_n \rightarrow 1$, we reobtain 
Theorem 1 in \cite{bm} and Theorem 3.2 in \cite{rouhanimoradi}. In addition if $\alpha_n=0$ for all $n\ge 0$ (or for all $n\ge N$) then (\ref{generalalg}) 
defines just a simple sequence (not an iterative method since $x_{n+1}$ is no longer dependent on $x_n$) which approximates $P_Fu$ and in this case Theorem 3.4 in \cite{rouhanimoradi} is reobtained as a simple particular case 
(with $u_n:= \lambda_n u + (1-\lambda_n)(y_0+e_n), \ n\ge 0$). 

\section{Main Result}
Since we want to show that the sequences generated by (\ref{generalalg}) are convergent, we begin this section with a preliminary result stating that a 
necessary condition is $F=A^{-1}0 \neq \emptyset$. 
\begin{lemma}\label{Theorem1}
Assume that $A:D(A)\subset H \to H$ is a maximal monotone operator, $\beta_n \rightarrow \infty$, $(u_n)_{n\ge 0}$ is a bounded sequence, 
$|\alpha_n|\le c \ \forall n\ge 0$ for some $c <1$, and $(\alpha_ne_n)$ is bounded. Then the sequence $(x_n)$ 
generated by (\ref{generalalg}) is bounded for all $x_0\in H$ (equivalently, for some $x_0\in H$) if and only if $F\neq \emptyset$.
\end{lemma}
{\it Proof} 
Assume that for some $x_0\in H$ the sequence $(x_n)$ generated by (\ref{generalalg}) is bounded. We have 
$$
Ax_n\ni z_n:=\frac{1}{\beta_{n-1}}(u_{n-1}+\alpha_{n-1}x_{n-1}+\alpha_{n-1}e_{n-1}-x_{n-1}) \rightarrow 0 \, . 
$$
Therefore taking the limit in the obvious inequality
$$
(v-x_n,w-z_n)\ge 0 \ \ \ \forall [v,w]\in G(A) \, , 
$$
we infer 
$$
(v-p,w-0)\ge 0 \ \ \ \forall [v,w]\in G(A) \, , 
$$
where $p$ is a weak cluster point of $(x_n)$. By the maximality of $A$ we obtain
$$
[p,0]\in G(A) \ \Rightarrow \ p\in D(A), \ 0\in Ap \, . 
$$ 
Conversely, assume $F\neq \emptyset$. Let $p\in F$ and $x_0\in H$ be arbitrary but fixed points. Since the resolvent operator is nonexpansive we have 
\begin{eqnarray}
\Vert x_{n+1}-p\Vert &=& \Vert J_{\beta_n}(u_n+\alpha_n(x_n+e_n))-J_{\beta_n}p\Vert \nonumber \\
&\le&\Vert J_{\beta_n}(u_n+\alpha_n(x_n+e_n))-J_{\beta_n}u_n\Vert + \Vert J_{\beta_n}u_n-J_{\beta_n}p\Vert \nonumber \\
&\le& |\alpha_n|(\Vert x_n\Vert + \Vert e_n\Vert)+ \Vert u_n-p\Vert \nonumber \\
&\le& c\Vert x_n\Vert +|\alpha_n|\cdot \Vert e_n\Vert + \Vert u_n\Vert + \Vert p\Vert \nonumber \, . 
\end{eqnarray}
Therefore
\begin{equation}\label{two}
\Vert x_{n+1}\Vert \le c\Vert x_n\Vert + c_1 \ \ \ \forall n\ge 0 \, ,  
\end{equation}
where $c_1$ is a positive constant. From (\ref{two}) we derive by induction 
\begin{eqnarray}
\Vert x_n\Vert &\le& \Vert x_0\Vert c^n + c_1\underbrace{(c^{n-1}+ c^{n-2} + \cdots +c+1)}_{=\frac{1-c^n}{1-c}} \nonumber \\
&\le& \Vert x_0\Vert c^n + \frac{c_1}{1-c} \ \ \forall n\ge 2 \, , \nonumber 
\end{eqnarray}
which shows that $(x_n)$ is bounded.
\qed

\bigskip 

Before stating our main theorem let us recall the following result which was proved independently by Bruck \cite{bruck} and Moro\c sanu \cite{moro}. 
\begin{lemma}\label{moro}
Let $A:D(A)\subset H \to H$ be a maximal monotone operator with $F=:A^{-1}0$ nonempty. Then for every $u\in H$, 
$(I+\lambda A)^{-1}u \rightarrow P_Fu$ as $\lambda \rightarrow \infty$, where $P_Fu$ denotes the metric projection of $u$ onto $F$. 
\end{lemma}
Now let us state the main result of this Note.
\begin{theorem}\label{main}
Assume $(H)$ (see the previous section). Then for all $x_0\in H$ the sequence $(x_n)$ generated by algorithm (\ref{generalalg}) converges 
strongly to $P_Fu$ (the metric projection of $u$ onto $F=A^{-1}0$). 
\end{theorem}
{\it Proof}
Let $x_0\in H$ be an arbitrary but fixed point. By Theorem \ref{Theorem1} the corresponding sequence $(x_n)$ generated by (\ref{generalalg})  
is bounded (since there exists a natural number $N$ such that $|\alpha_n|\le c<1$ for $n\ge N$ so Lemma  \ref{Theorem1} is applicable with $x_0:=x_N$). 
Thus we have 
\begin{eqnarray}
\Vert x_{n+1}-P_Fu\Vert &=& \Vert J_{\beta_n}(u_n+\alpha_n(x_n+e_n))-J_{\beta_n}u\Vert + \Vert J_{\beta_n}u-P_Fu\Vert \nonumber \\
&\le&\Vert (u_n-u)+\alpha_n(x_n + e_n)\Vert + \Vert J_{\beta_n}u_n-P_Fu\Vert \nonumber \\
&\le& \Vert u_n-u\Vert +|\alpha_n|\cdot \Vert x_n\Vert + |\alpha_n|\cdot \Vert e_n\Vert + \Vert J_{\beta_n}u-P_Fu\Vert \, . \nonumber 
\end{eqnarray}
So by $(H)$ and Lemma \ref{moro} we conclude that $\Vert x_n - P_Fu\Vert \rightarrow 0$ as $n\rightarrow \infty$. 
\qed

\section{Concluding Comments}

1. If $A=\partial \phi$ where $\phi:H\to (-\infty , \infty]$ is a proper, convex, lower semicontinuous function the algorithm (\ref{generalalg}) serves 
as a method for appoximating minimizers of $\phi$ (assuming that the set of minimizers of $\phi$ is nonempty), since 
in this case $p\in A^{-1}0$ if and only if $p$ is a minimizer of $\phi$. 

\bigskip

2. The error sequence $(e_n)$ is allowed to be bounded as usual in numerical analysis, or even unbounded with $\alpha_n\Vert e_n\Vert \rightarrow 0$. 

\bigskip

3. Theorem \ref{main} is a generalization of both Theorem 3.2 in \cite{rouhanimoradi} and Theorem 1 in \cite{bm}. 

\bigskip

If $\alpha_n=0$ for all $n\ge 0$ (or for all $\alpha_n\ge N$) then (\ref{generalalg}) defines just a simple sequence, not a real iterative method, since 
$x_{n+1}$ does not depend on $x_n$. In this case we have 
\begin{equation}\label{rou}
x_{n+1}=J_{\beta_n}u_n \ \ \ \forall n\ge 0 \, , 
\end{equation}
or for all $n\ge N$. In fact, according to Lemma \ref{moro}, we have 
\begin{eqnarray}
\Vert x_{n+1} -P_Fu\Vert &\le & \Vert J_{\beta_n}u_n - J_{\beta_n}u\Vert + \Vert J_{\beta_n}u - P_Fu\Vert \nonumber \\
&\le& \Vert u_n-u\Vert + \Vert J_{\beta_n}u - P_Fu\Vert \rightarrow 0 \, . \nonumber 
\end{eqnarray}
Note that the second algorithm introduced and studied in \cite[p. 228]{rouhanimoradi} is in fact a sequence of the form (\ref{rou}) with 
$u_n=\lambda_nu + (1-\lambda_n)(y_0+e_n)$, $n\ge 0$. 

5. Let us explain how the algorithn (\ref{generalalg}) works when performing simulations. Assume conditions $(H)$ are fulfilled. In addition, for 
the sake of simplicity, $A$ is assumed to be single-valued. For a given $x_0\in H$ we compute $x_1$ by solving for $x$ the equation
$$
(I+\beta_0A)x = u_0 + \alpha_0x_0 
$$
and get $x_1 + e_1$ instead of the exact solution $x=x_1$. We do not have any error for $x_0$ (i.e., $e_0=0$) but we have a computational 
error $e_1$ for $x_1$. Next we solve for $x$ the equation
$$
(I+\beta_1A)x = u_1 + \alpha_1(x_1 + e_1)
$$
and get $x_2 + e_2$ instead of the exact solution 
$$
x_2=(I+\beta_1A)^{-1}(u_1+\alpha_1(x_1+e_1)) \, ,
$$
and so on. Thus using the computer we obtain $z_n=x_n+e_n$ satisfying 
$$
z_{n+1}=(I+\beta_nA)^{-1}(u_n + \alpha_nz_n) + e_{n+1} \ \ \ \forall n\ge 0\, , 
$$
where $z_0=0$. If $\Vert e_n\Vert \le \varepsilon$ for all $n\ge 0$, then $\Vert x_n-z_n\Vert \le \varepsilon$ for all $n\ge 0$, i.e., 
$z_n$ approximates $P_Fu$ for $n$ large enough. 

\section{Simulations}

Intuitively, in order to achieve fast convergence of the sequence $(x_n)$ generated by algorithm (\ref{generalalg}) to 
$P_Fu$ we need to choose a point $u$ as close as possible to $F=A^{-1}0$ and sequences $(\beta_n)$ and $(\alpha_n)$ that 
converge fastly to $\infty$ and to $0$, respectively.

etc., etc., ... to be continued. 



\begin{thebibliography}{}

\bibitem{brezis} Brezis, H.: Op\'{e}rateurs Maximaux Monotones et Semi-Groupes de Contractions dans les Espaces de Hilbert, North Holland Math. 
Studies, Vol. 5. North Holland, Amsterdam (1973) 

\bibitem{morosanu} Moro\c{s}anu, G.: Nonlinear Evolution Equations and Applications. Reidel, Dordrecht (1988)

\bibitem{martinet} Martinet, B.: R\'{e}gularisation d'in\'{e}quations variationnelles par approximations succesives. Rev. Fran\c{c}aise Informat. 
Recherche Op\'{e}rationnelle {\bf 4} (Ser. R-3), 154-158 (1970)

\bibitem{rockafellar} Rockafellar, R.T.: Monotone operators and the proximal point algorithm. SIAM J. Control Optim. {\bf 14}, 877-898 (1976)

\bibitem{guler} G\"{u}ler, O.: On the convergence of the proximal point algorithm for convex minimization. SIAM J. Control Optim. {\bf 29}, 403-419 (1991)

\bibitem{lehmud} Lehdili, N., Moudafi, A.: Combining the proximal algorithm and Tikhonov regularization. Optimization {\bf 37}, 239-252 (1996)

\bibitem{xu} Xu, H.K.: A regularization method for the proximal point algorithm. J. Global Optim. {\bf 36}, 115-125 (2006)

\bibitem{wangcui} Wang, F., Cui, H.: On the contraction-proximal point algorithms with multi-parameters. J. Global Optim. {\bf 54}, 485-491 (2012)

\bibitem{bm} Boikanyo, O.A., Moro\c{s}anu, G.: Strong convergence of a proximal point algorithm with bounded error sequence. Optim. Lett. {\bf 7}, 
415-420 (2013)

\bibitem{rouhanimoradi} Djafari Rouhani, B., Moradi, S.: Strong convergence of two proximal point algorithms with possible unbounded error sequences.  
J. Optim. Theory Appl. {\bf 172}, 222-235 (2017)

\bibitem{bruck} Bruck, R.E. Jr.: A strongly convergent iterative solution of $0\in U(x)$ for a maximal monotone operator $U$ in Hilbert space. 
J. Math. Anal. Appl. {\bf 48}, 114-126 (1974) 

\bibitem{moro} Moro\c{s}anu, G.: Asymptotic behaviour of resolvent for a monotone mapping in a Hilbert space. Atti Accad. Naz. Lincei Rend. 
Cl. Sci. Fis. Mat. Natur. {\bf 61}, 565-570 (1977)

\end{thebibliography}
\end{document}